\documentclass[11pt]{article}
\usepackage{amsfonts,amscd,amssymb, amsmath}
\usepackage[all]{xy}
\newtheorem{definition}{Definition}[section]
\newtheorem{lemma}[definition]{Lemma}
\newtheorem{proposition}[definition]{Proposition}

\newtheorem{remark}[definition]{Remark}

\newtheorem{theorem}[definition]{Theorem}

\def\rawo\lonra{\longrightarrow}

\def\ot{\otimes}

\allowdisplaybreaks[4]

\newenvironment{proof}{{\it Proof.}}{\hfill $ \square $ \vskip 4mm}
\begin{document}
\title{More examples of pseudosymmetric braided categories
\thanks{Research partially supported by the  
CNCSIS project ''Hopf algebras, cyclic homology and monoidal categories'', 
contract nr. 560/2009, CNCSIS code $ID_{-}69$.}}
\author{Florin Panaite\\
Institute of Mathematics of the 
Romanian Academy\\ 
PO-Box 1-764, RO-014700 Bucharest, Romania\\
e-mail: Florin.Panaite@imar.ro
\and 
Mihai D. Staic\thanks{Institute of Mathematics of the  Romanian Academy, 
PO-Box 1-764, RO-014700 Bucharest, Romania.}\\
Department of Mathematics and Statistics, BGSU\\
Bowling Green, OH 43403, USA\\
e-mail: mstaic@bgsu.edu} 
\date{}
\maketitle
%%%%%%%%%%%%%%%%%%%%%%%%%%%%%%%
\begin{abstract}
 We study some examples of braided categories and quasitriangular Hopf algebras and 
decide which of them is pseudosymmetric, respectively pseudotriangular. We show also 
that there exists a universal pseudosymmetric braided category.
\end{abstract}
%%%%%%%%%%%%%%%%%%%%%%%%%%%%%%%%%%
%%%%%%%%%%%%%%%%%%%%%%%%%%%%%%%%%%%
\section*{Introduction}
%%%%%%%%%%%%%%%%%%%%%%%%%%%%%%%%%%%%
${\;\;\;\;}$
Braided categories have been  introduced by Joyal and Street in \cite{js} as natural 
generalizations of symmetric categories. Roughly speaking, a braided category is a category  
that has a tensor product with a nice commutation rule. More precisely, for every two 
objects $U$ and $V$ we have an isomorphism $c_{U,V}:U\otimes V\to V\otimes U$ that 
satisfies certain conditions. These conditions are chosen in such a way that for every 
object $V$ in the category there exists a natural way to construct  a representation 
for the braid group $B_n$ on  $V^{\otimes n}$, therefore the name {\it braided } categories. 
If we impose the extra condition $c_{V,U}c_{U,V}=id_{U\otimes V}$ for all objects $U, V$ in the category, we recover 
the definition of symmetric categories. It is well known that symmetric categories 
can be used to construct representations for the symmetric group $\Sigma_n$.   

Pseudosymmetric categories  are a special class of braided categories and have been 
introduced in \cite{psvo}. The motivation was the study of certain categorical structures 
called twines, strong twines and pure-braided structures (introduced in \cite{brug}, \cite{psv} 
and  \cite{doru}). A braiding on a strict monoidal category is called pseudosymmetric if it 
satisfies a sort of modified braid relation; any symmetric braiding is pseudosymmetric.   
One of the most intriguing results obtained in \cite{psvo} was that the category of 
Yetter-Drinfeld modules over a Hopf algebra $H$ is pseudosymmetric if and only if $H$ is commutative and cocommutative.   We proved in \cite{ps} that pseudosymmetric categories 
can be used to construct representations for the group $PS_n=\frac{B_n}{[P_n,P_n]}$, 
the quotient of the braid group by the commutator subgroup of the pure braid group.  
There exists also a Hopf algebraic analogue of pseudosymmetric braidings: a quasitriangular 
structure on a Hopf algebra is called pseudotriangular if it satisfies a sort of modified
quantum Yang-Baxter equation. 

In this paper we tie some lose ends from \cite{psvo} and \cite{ps}. We study  more examples 
of braided categories and quasitriangular Hopf algebras and decide when they are pseudosymmetric, respectively pseudotriangular. Namely, we prove that the canonical braiding  
of the category ${\cal LR}(H)$  of Yetter-Drinfeld-Long bimodules over a Hopf algebra $H$ 
(introduced in \cite{pvo}) is pseudosymmetric if and only if 
$H$ is commutative and cocommutative. We show that any quasitriangular structure on the 
$4\nu $-dimensional Radford's Hopf algebra $H_{\nu }$ (introduced in \cite{radford}) 
is pseudotriangular.  We analyze the positive quasitriangular structures $R(\xi , \eta )$ on a 
Hopf algebra with positive bases $H(G; G_+, G_-)$ (as defined in \cite{LYZ1}, \cite{LYZ2}), 
where $\xi , \eta $ are group homomorphisms from $G_+$ to $G_-$, and we present a list of 
necessary and sufficient conditions for $R(\xi , \eta )$ to be pseudotriangular. If 
$R(\xi , \eta )$ is normal (i.e. if $\xi $ is trivial) these conditions reduce to the single relation 
$\eta (uv)=\eta (vu)$ for all $u, v\in G_+$. 

In the last section we recall the pseudosymmetric braided category ${\cal PS}$ introduced in \cite{ps} and we show that it is a universal pseudosymmetric category. 
More precisely, we prove that it  
satisfies two  universality properties similar to the ones satisfied by the 
universal braid category ${\cal B}$ (see \cite{k}).   
%%%%%%%%%%%%%%%%%%%%%%%%%%%%
\section{Preliminaries}\label{sec1}
%%%%%%%%%%%%%%%%%%%%%%%%%%%%%%%
\setcounter{equation}{0}
%%%%%%%%%%%%%%%%%%%%%%%%
${\;\;\;\;}$
We work over a base field $k$. 
All algebras, linear 
spaces, etc, will be over $k$; unadorned $\ot$ means $\ot _k$. 
For a Hopf algebra $H$ with comultiplication $\Delta$ we  
denote $\Delta (h)=h_1\ot h_2$, for $h\in H$. 
For terminology concerning Hopf algebras and  
monoidal categories we refer to \cite{k}. 
\begin{definition} (\cite{psvo}) 
Let ${\mathcal C}$ be a strict monoidal category and $c$ a braiding on  
${\mathcal C}$. We say that $c$ is {\bf pseudosymmetric} if the following   
condition holds, for all $X, Y, Z\in {\mathcal C}$:
\begin{eqnarray*}
(c_{Y, Z}\otimes id_X)(id_Y\otimes c_{Z, X}^{-1})
(c_{X, Y}\otimes id_Z)
=(id_Z\otimes c_{X, Y})(c_{Z, X}^{-1}\otimes id_Y)
(id_X\otimes c_{Y, Z}). \nonumber
\end{eqnarray*}
In this case we say that ${\mathcal C}$ is   
a {\bf pseudosymmetric braided category}.  
\end{definition}
\begin{proposition} (\cite{psvo})
Let ${\mathcal C}$ be a strict monoidal category and $c$ a braiding on  
${\mathcal C}$. Then $c$ is pseudosymmetric if and only if the family 
$T_{X, Y}:=c_{Y, X}c_{X, Y}:X\otimes Y\rightarrow X\otimes Y$ satisfies the condition  
$(T_{X, Y}\otimes id_Z)(id_X\otimes T_{Y, Z})=
(id_X\otimes T_{Y, Z})(T_{X, Y}\otimes id_Z)$ for all $X, Y, Z\in {\cal C}$. 
\end{proposition}
\begin{definition} (\cite{psvo}) 
Let $H$ be a Hopf algebra and $R\in H\otimes H$ a quasitriangular 
structure. Then $R$ is called {\bf pseudotriangular} if 
$R_{12}R^{-1}_{31}R_{23}=R_{23}R^{-1}_{31}R_{12}$.
\end{definition}
\begin{proposition} (\cite{psvo})
Let $H$ be a Hopf algebra and let R be a quasitriangular structure on $H$. Then $R$ is 
pseudotriangular 
if and only if the element $F=R_{21}R\in H\otimes H$ satisfies the relation
$F_{12}F_{23}=F_{23}F_{12}$.
\end{proposition}

%%%%%%%%%%%%%%%%%%%%%%%%%%%%%%%%%%%%
\section{Yetter-Drinfeld-Long bimodules} 
%%%%%%%%%%%%%%%%%%%%%%%%%%%%%%%%
\setcounter{equation}{0}
%%%%%%%%%%%%%%%%%%%%%%%%%%%%%%
${\;\;\;\;}$For a braided monoidal category ${\mathcal C}$ with braiding $c$, 
let ${\mathcal C}^{in}$ be equal to ${\mathcal C}$ as a monoidal category, 
with the mirror-reversed braiding $\tilde{c}_{M, N}:=c_{N, M}^{-1}$, 
for all objects $M, N\in {\mathcal C}$. Directly from the definition of a 
pseudosymmetric braiding, we immediately obtain: 
\begin{proposition} \label{inversepseudo}
Let ${\mathcal C}$ be a strict braided monoidal category. Then ${\mathcal C}$ 
is pseudosymmetric if and only if ${\mathcal C}^{in}$ 
is pseudosymmetric. 
\end{proposition}

Let $H$ be a Hopf algebra with bijective antipode $S$. 
Consider the category $_H{\cal YD}^H$ of 
left-right Yetter-Drinfeld modules over $H$, whose objects are vector 
spaces $M$ that are left $H$-modules (denote the action by 
$h\otimes m\mapsto h\cdot m$) and right $H$-comodules (denote the coaction 
by $m\mapsto m_{(0)}\otimes m_{(1)}\in M\otimes H$) satisfying the 
compatibility condition 
\begin{eqnarray*}
&&(h\cdot m)_{(0)}\otimes (h\cdot m)_{(1)}=h_2\cdot m_{(0)}\otimes 
h_3m_{(1)}S^{-1}(h_1), \;\;\;\forall \;\;h\in H, \;m\in M. 
\end{eqnarray*}
It is a monoidal category, with tensor product given by 
\begin{eqnarray*}
&&h\cdot (m\otimes n)=h_1\cdot m\otimes h_2\cdot n, \;\;\;
(m\otimes n)_{(0)}\otimes (m\otimes n)_{(1)}=m_{(0)}\otimes n_{(0)}
\otimes n_{(1)}m_{(1)}.
\end{eqnarray*}
Moreover, it has a (canonical) braiding given by 
\begin{eqnarray*}
&&c_{M, N}:M\otimes N\rightarrow N\otimes M, \;\;\;c_{M, N}(m\otimes n)=
n_{(0)}\otimes n_{(1)}\cdot m, \\
&&c_{M, N}^{-1}:N\otimes M\rightarrow M\otimes N, \;\;\;
c_{M, N}^{-1}(n\otimes m)=S(n_{(1)})\cdot m\otimes n_{(0)}. 
\end{eqnarray*}

Consider also the category $_H^H{\cal YD}$ of 
left-left Yetter-Drinfeld modules over $H$, whose objects are vector 
spaces $M$ that are left $H$-modules (denote the action by 
$h\otimes m\mapsto h\cdot m$) and left $H$-comodules (denote the coaction 
by $m\mapsto m^{(-1)}\otimes m^{(0)}\in H\otimes M$) with
compatibility condition 
\begin{eqnarray*}
&&(h_1\cdot m)^{(-1)}h_2\otimes (h_1\cdot m)^{(0)}=h_1m^{(-1)}\otimes 
h_2\cdot m^{(0)}, \;\;\;\forall \;\;h\in H, \;m\in M. 
\end{eqnarray*}
It is a monoidal category, with tensor product given by 
\begin{eqnarray*}
&&h\cdot (m\otimes n)=h_1\cdot m\otimes h_2\cdot n, \;\;\;
(m\otimes n)^{(-1)}\otimes (m\otimes n)^{(0)}=m^{(-1)}n^{(-1)}
\otimes m^{(0)}\otimes n^{(0)}.
\end{eqnarray*}
Moreover, it has a (canonical) braiding given by 
\begin{eqnarray*}
&&c_{M, N}:M\otimes N\rightarrow N\otimes M, \;\;\;c_{M, N}(m\otimes n)=
m^{(-1)}\cdot n\otimes m^{(0)}, \\
&&c_{M, N}^{-1}:N\otimes M\rightarrow M\otimes N, \;\;\;
c_{M, N}^{-1}(n\otimes m)=m^{(0)}\otimes S^{-1}(m^{(-1)})\cdot n. 
\end{eqnarray*}
\begin{proposition} (\cite{bcp}) \label{isoyd}
For the categories $_H{\cal YD}^H$ and $_H^H{\cal YD}$ with 
braidings as above, we have an isomorphism of braided monoidal categories 
$(_H{\cal YD}^H)^{in}\simeq \;$$_H^H{\cal YD}$.
\end{proposition}
\begin{proposition} (\cite{psvo}) \label{charact}
The canonical braiding of $_H{\cal YD}^H$ is pseudosymmetric if and only if 
$H$ is commutative and cocommutative. 
\end{proposition}

As a consequence of Propositions \ref{inversepseudo}, \ref{isoyd} and 
\ref{charact}, we obtain:
\begin{proposition}\label{ydleft}
The canonical braiding of $_H^H{\cal YD}$ is pseudosymmetric if and only if 
$H$ is commutative and cocommutative. 
\end{proposition}

We recall now the braided monoidal category ${\cal LR}(H)$ 
defined in \cite{pvo}.
The objects of 
${\cal LR}(H)$ are vector spaces $M$ endowed with $H$-bimodule and 
$H$-bicomodule structures (denoted by $h\otimes m\mapsto h\cdot m$, 
$m\otimes h\mapsto m\cdot h$, $m\mapsto m^{(-1)}\otimes m^{(0)}$, 
$m\mapsto m^{<0>}\otimes m^{<1>}$, for all $h\in H$, $m\in M$), such that 
$M$ is a left-left Yetter-Drinfeld module, a left-right Long module, 
a right-right Yetter-Drinfeld module and a right-left Long module, i.e.  
(for all $h\in H$, $m\in M$):
\begin{eqnarray}
&&(h_1\cdot m)^{(-1)}h_2\otimes (h_1\cdot m)^{(0)}=
h_1m^{(-1)}\otimes h_2\cdot m^{(0)}, \label{caty1} \\
&&(h\cdot m)^{<0>}\otimes (h\cdot m)^{<1>}=h\cdot m^{<0>}\otimes 
m^{<1>}, \label{caty2} \\
&&(m\cdot h_2)^{<0>}\otimes h_1(m\cdot h_2)^{<1>}=
m^{<0>}\cdot h_1\otimes m^{<1>}h_2, \label{caty3} \\
&&(m\cdot h)^{(-1)}\otimes (m\cdot h)^{(0)}=m^{(-1)}\otimes 
m^{(0)}\cdot h. \label{caty4}
\end{eqnarray}
Morphisms in ${\cal LR}(H)$ are 
$H$-bilinear $H$-bicolinear maps. ${\cal LR}(H)$ is
a strict monoidal category, 
with unit $k$ endowed with usual $H$-bimodule and $H$-bicomodule structures, 
and tensor product given by: if $M, N\in {\cal LR}(H)$ then 
$M\otimes N\in {\cal LR}(H)$ as follows (for all $m\in M$, $n\in N$, 
$h\in H$): 
\begin{eqnarray*}
&&h\cdot (m\otimes n)=h_1\cdot m\otimes h_2\cdot n, \;\;\;
(m\otimes n)\cdot h=m\cdot h_1\otimes n\cdot h_2, \\
&&(m\otimes n)^{(-1)}\otimes (m\otimes n)^{(0)}=m^{(-1)}n^{(-1)}\otimes 
(m^{(0)}\otimes n^{(0)}), \\
&&(m\otimes n)^{<0>}\otimes (m\otimes n)^{<1>}=(m^{<0>}\otimes n^{<0>})
\otimes m^{<1>}n^{<1>}.
\end{eqnarray*}
Moreover, ${\cal LR}(H)$ has a (canonical) braiding defined, for 
$M, N\in {\cal LR}(H)$, $m\in M$, $n\in N$, by 
\begin{eqnarray*}
&&c_{M, N}:M\otimes N\rightarrow N\otimes M, \;\;\;
c_{M, N}(m\otimes n)=m^{(-1)}\cdot n^{<0>}\otimes m^{(0)}\cdot n^{<1>}, \\
&&c_{M, N}^{-1}:N\otimes M\rightarrow M\otimes N, \;\;\;
c_{M, N}^{-1}(n\otimes m)=m^{(0)}\cdot S^{-1}(n^{<1>})\otimes 
S^{-1}(m^{(-1)})\cdot n^{<0>}. 
\end{eqnarray*}
\begin{proposition}
The canonical braiding of ${\cal LR}(H)$  is pseudosymmetric if and 
only if 
$H$ is commutative and cocommutative. 
\end{proposition}
\begin{proof}
Assume that the canonical braiding of ${\cal LR}(H)$  is pseudosymmetric. 
As noted in \cite{pvo}, $_H^H{\cal YD}$ with its canonical braiding is a 
braided subcategory of ${\cal LR}(H)$, so the canonical braiding of 
$_H^H{\cal YD}$ is pseudosymmetric; by Proposition \ref{ydleft} 
it follows that $H$ is commutative and cocommutative.\\
Conversely, assume that $H$ is commutative and cocommutative. Then one can see 
that the two Yetter-Drinfeld conditions appearing in the definition of 
${\cal LR}(H)$ 
become Long conditions, that is (\ref{caty1}) and (\ref{caty3}) become 
respectively
\begin{eqnarray}
&&(h\cdot m)^{(-1)}\otimes (h\cdot m)^{(0)}=
m^{(-1)}\otimes h\cdot m^{(0)}, \label{longnew1} \\
&&(m\cdot h)^{<0>}\otimes (m\cdot h)^{<1>}=
m^{<0>}\cdot h\otimes m^{<1>}. \label{longnew3}
\end{eqnarray}
Let now $X, Y, Z\in {\cal LR}(H)$; we compute, for $x\in X$, 
$y\in Y$, $z\in Z$:\\[2mm]
${\;\;\;\;\;}$$(c_{Y, Z}\otimes id_X)(id_Y\otimes c_{Z, X}^{-1})
(c_{X, Y}\otimes id_Z)(x\ot y\ot z)$
\begin{eqnarray*}
&=&(c_{Y, Z}\otimes id_X)(id_Y\otimes c_{Z, X}^{-1})
(x^{(-1)}\cdot y^{<0>}\ot x^{(0)}\cdot y^{<1>}\ot z)\\
&=&(c_{Y, Z}\otimes id_X)(x^{(-1)}\cdot y^{<0>}\ot 
z^{(0)}\cdot S^{-1}((x^{(0)}\cdot y^{<1>})^{<1>})\\
&&\ot 
S^{-1}(z^{(-1)})\cdot (x^{(0)}\cdot y^{<1>})^{<0>})\\
&\overset{(\ref{longnew3})}{=}&
(c_{Y, Z}\otimes id_X)(x^{(-1)}\cdot y^{<0>}\ot 
z^{(0)}\cdot S^{-1}(x^{(0)<1>})\ot 
S^{-1}(z^{(-1)})\cdot x^{(0)<0>}\cdot y^{<1>})\\
&=&(x^{(-1)}\cdot y^{<0>})^{(-1)}\cdot [z^{(0)}\cdot 
S^{-1}(x^{(0)<1>})]^{<0>} \\
&&\ot (x^{(-1)}\cdot y^{<0>})^{(0)}\cdot [z^{(0)}\cdot 
S^{-1}(x^{(0)<1>})]^{<1>} \\
&&\ot S^{-1}(z^{(-1)})\cdot x^{(0)<0>}\cdot y^{<1>}\\
&\overset{(\ref{longnew1}, \ref{longnew3})}{=}&
y^{<0>(-1)}\cdot z^{(0)<0>}\cdot S^{-1}(x^{(0)<1>})\ot 
x^{(-1)}\cdot y^{<0>(0)}\cdot z^{(0)<1>}\\
&&\ot 
S^{-1}(z^{(-1)})\cdot x^{(0)<0>}\cdot y^{<1>}, 
\end{eqnarray*}
${\;\;\;\;\;}$$(id_Z\otimes c_{X, Y})(c_{Z, X}^{-1}\otimes id_Y)
(id_X\otimes c_{Y, Z})(x\ot y\ot z)$
\begin{eqnarray*}
&=&(id_Z\otimes c_{X, Y})
(c_{Z, X}^{-1}\otimes id_Y)(x\ot 
y^{(-1)}\cdot z^{<0>}\ot y^{(0)}\cdot z^{<1>})\\
&=&(id_Z\otimes c_{X, Y})([y^{(-1)}\cdot z^{<0>}]^{(0)}
\cdot S^{-1}(x^{<1>})\ot S^{-1}([y^{(-1)}\cdot z^{<0>}]^{(-1)})
\cdot x^{<0>}\\
&&\ot y^{(0)}\cdot z^{<1>})\\
&\overset{(\ref{longnew1})}{=}&
(id_Z\otimes c_{X, Y})(y^{(-1)}\cdot z^{<0>(0)}
\cdot S^{-1}(x^{<1>})\ot S^{-1}(z^{<0>(-1)})
\cdot x^{<0>}\ot y^{(0)}\cdot z^{<1>})\\
&=&y^{(-1)}\cdot z^{<0>(0)}\cdot S^{-1}(x^{<1>})\ot 
[S^{-1}(z^{<0>(-1)})\cdot x^{<0>}]^{(-1)}\cdot 
[y^{(0)}\cdot z^{<1>}]^{<0>}\\
&&\ot [S^{-1}(z^{<0>(-1)})\cdot x^{<0>}]^{(0)}\cdot 
[y^{(0)}\cdot z^{<1>}]^{<1>}\\
&\overset{(\ref{longnew1}, \ref{longnew3})}{=}&
y^{(-1)}\cdot z^{<0>(0)}\cdot S^{-1}(x^{<1>})\ot 
x^{<0>(-1)}\cdot y^{(0)<0>}\cdot z^{<1>}\\
&&\ot 
S^{-1}(z^{<0>(-1)})\cdot x^{<0>(0)}\cdot 
y^{(0)<1>}, 
\end{eqnarray*}
and the two terms are equal because of the bicomodule 
condition for $X$, $Y$ and $Z$.
\end{proof}
%%%%%%%%%%%%%%%%%%%%%%%%%%%%%%%%%%%%
\section{Radford's Hopf algebras $H_{\nu }$} 
%%%%%%%%%%%%%%%%%%%%%%%%%%%%%%%%
\setcounter{equation}{0}
%%%%%%%%%%%%%%%%%%%%%%%%%%%%%%
${\;\;\;\;}$
Let $\nu $ be an odd natural number and assume that the base field $k$
contains a primitive $2\nu ^{th}$ root of unity $\omega $ and 
$2\nu $ is invertible in $k$. We consider a certain family of Hopf 
algebras, which are exactly the quasitriangular ones from the 
larger family constructed by Radford in \cite{radford}. Namely, 
using notation as in \cite{cc}, we denote by $H_{\nu }$ the 
Hopf algebra over $k$ generated by two elements $g$ and $x$ 
such that 
\begin{eqnarray*}
&&g^{2\nu }=1, \;\;\;gx+xg=0, \;\;\;x^2=0, 
\end{eqnarray*}
with coproduct $\Delta (g)=g\ot g$ and $\Delta (x)=x\ot g^{\nu }+1\ot x$, and 
antipode $S(g)=g^{-1}$ and $S(x)=g^{\nu }x$. Note that $H_1$ is exactly 
Sweedler's 4-dimensional Hopf algebra, and in general 
$H_{\nu }$ is $4\nu $-dimensional, a linear basis in $H_{\nu }$ being the set 
$\{g^lx^m/0\leq l<2\nu , \;0\leq m\leq 1\}$. 

The quasitriangular structures of $H_{\nu }$ have been determined in \cite{radford}; 
they are parametrized by pairs $(s, \beta )$, where $\beta \in k$ and 
$s$ is an odd number with $1\leq s<2\nu $. Moreover, if we denote by 
$R_{s, \beta }$ the quasitriangular structure corresponding to 
$(s, \beta )$, then we have 
\begin{eqnarray*}
&&R_{s, \beta }=\frac{1}{2\nu }(\sum _{i, l=0}^{2\nu -1}
\omega ^{-il}g^i\ot g^{sl})+\frac{\beta }{2\nu }(\sum _{i, l=0}^{2\nu -1}
\omega ^{-il}g^ix\ot g^{sl+\nu }x).
\end{eqnarray*}
It was also proved in \cite{radford} that $R_{s, \beta }$ is triangular if and only if 
$s=\nu $.  

Following \cite{radford}, we introduce an alternative description of $R_{s, \beta }$, 
more appropriate for our purpose. For every natural number 
$0\leq l\leq 2\nu -1$, we define
\begin{eqnarray*} 
&&e_l=\frac{1}{2\nu }\displaystyle\sum _{i=0}^{2\nu -1}
\omega ^{-il}g^i, 
\end{eqnarray*}
regarded as an element in the group algebra of the 
cyclic group of order $2\nu $ generated by the element $g$ (which in turn 
may be regarded as a Hopf subalgebra of $H_{\nu }$ in the obvious way). Then, 
by \cite{radford}, the following relations hold: 
\begin{eqnarray*}
&&1=e_0+e_1+...+e_{2\nu -1}, \label{idem1}\\
&&e_ie_j=\delta _{ij}e_i, \label{idem2}\\
&&g^ie_j=\omega ^{ij}e_j, \label{idem3} 
\end{eqnarray*}
for all $0\leq i, j\leq 2\nu -1$. Also, a straightforward computation shows that 
we have 
\begin{eqnarray*}
&&\displaystyle\sum _{i=0}^{2\nu -1}(-1)^ie_i=g^{\nu }. \label{idem4}
\end{eqnarray*}
Note also that, since $\omega $ is a primitive $2\nu ^{th}$ root of unity,  we 
have 
\begin{eqnarray*}
\omega ^{\nu }=-1. \label{radprim}
\end{eqnarray*}
With this notation, the quasitriangular structure $R_{s, \beta }$ may be expressed 
(cf. \cite{radford}) as 
\begin{eqnarray*}
&&R_{s, \beta }=\sum _{l=0}^{2\nu -1}
e_l\ot g^{sl}+\beta (\sum _{l=0}^{2\nu -1}
e_lx\ot g^{sl+\nu }x).
\end{eqnarray*}

We are interested to see for what $s, \beta $ is $R_{s, \beta }$ pseudotriangular. 
We note first that for $\beta =0$, $R_{s, 0}$ is actually a quasitriangular 
structure on the group algebra of the cyclic group of order $2\nu $, which is a 
commutative Hopf algebra, so $R_{s, 0}$ is pseudotriangular. 

Consider now $R_{s, \beta }$ an arbitrary quasitriangular structure on $H_{\nu }$. 
We need to compute first $(R_{s, \beta })_{21}R_{s, \beta }$. By using the 
defining relations $x^2=0$ and $gx+xg=0$, the properties of the elements 
$e_l$ listed above and the fact that $s$ and $\nu $ are odd numbers, a 
straightforward computation yields:
\begin{eqnarray*}
(R_{s, \beta })_{21}R_{s, \beta }&=&\sum _{l, t=0}^{2\nu -1}
\omega ^{2slt}e_l\ot e_t+\beta (\sum _{l, t=0}^{2\nu -1}
\omega ^{2slt+\nu t}e_lx\ot e_tx)\\
&&-\beta (\sum _{l, t=0}^{2\nu -1}
(-1)^{l+t}\omega ^{2slt+\nu l}xe_l\ot e_tx).
\end{eqnarray*}
Let us denote this element by $T$. We need to compare $T_{12}T_{23}$ and 
$T_{23}T_{12}$, so we first compute them, using repeatedly the 
defining relations of $H_{\nu }$ and the properties of the elements 
$e_l$:
\begin{eqnarray*}
T_{12}T_{23}&=&\sum _{l, t, j=0}^{2\nu -1}
\omega ^{2slt+2stj}e_l\ot e_t\ot e_j+\beta 
(\sum _{l, t, i, j=0}^{2\nu -1}
\omega ^{2slt+2sij+\nu j}e_l\ot e_te_ix\ot e_jx \\
&&-\sum _{l, t, i, j=0}^{2\nu -1}(-1)^{i+j}
\omega ^{2slt+2sij+\nu i}e_l\ot e_txe_i\ot e_jx
+\sum _{l, t, i, j=0}^{2\nu -1}
\omega ^{2slt+2sij+\nu t}e_lx\ot e_txe_i\ot e_j\\
&&-\sum _{l, t, i, j=0}^{2\nu -1}(-1)^{t+l}
\omega ^{2slt+2sij+\nu l}xe_l\ot e_txe_i\ot e_j)\\
&=&\sum _{l, t, j=0}^{2\nu -1}
\omega ^{2slt+2stj}e_l\ot e_t\ot e_j+\beta 
(\sum _{l, t, j=0}^{2\nu -1}(-1)^j
\omega ^{2slt+2stj}e_l\ot e_tx\ot e_jx \\
&&-\sum _{l, t, i, j=0}^{2\nu -1}(-1)^{j}
\omega ^{2slt+2sij}e_l\ot e_txe_i\ot e_jx
+\sum _{l, t, i, j=0}^{2\nu -1}(-1)^t
\omega ^{2slt+2sij}e_lx\ot e_txe_i\ot e_j\\
&&-\sum _{l, t, i, j=0}^{2\nu -1}(-1)^{t}
\omega ^{2slt+2sij}xe_l\ot e_txe_i\ot e_j)\\
&=&\sum _{l, t, j=0}^{2\nu -1}
\omega ^{2slt+2stj}e_l\ot e_t\ot e_j+\beta 
(\sum _{l, t, j=0}^{2\nu -1}(-1)^j
\omega ^{2slt}e_l\ot e_tx\ot g^{2st}e_jx \\
&&-\sum _{l, t, i, j=0}^{2\nu -1}(-1)^{j}
g^{2st}e_l\ot e_txe_i\ot g^{2si}e_jx
+\sum _{l, t, i, j=0}^{2\nu -1}(-1)^t
\omega ^{2sij}e_lx\ot g^{2sl}e_txe_i\ot e_j\\
&&-\sum _{l, t, i, j=0}^{2\nu -1}(-1)^{t}
\omega ^{2sij}xe_l\ot g^{2sl}e_txe_i\ot e_j)\\
&=&\sum _{l, t, j=0}^{2\nu -1}
\omega ^{2slt+2stj}e_l\ot e_t\ot e_j+\beta 
(\sum _{l, t=0}^{2\nu -1}
\omega ^{2slt}e_l\ot e_tx\ot g^{2st+\nu }x \\
&&-\sum _{t, i=0}^{2\nu -1}
g^{2st}\ot e_txe_i\ot g^{2si+\nu }x
+\sum _{l, i, j=0}^{2\nu -1}
\omega ^{2sij}e_lx\ot g^{2sl+\nu }xe_i\ot e_j\\
&&-\sum _{l, i, j=0}^{2\nu -1}
\omega ^{2sij}xe_l\ot g^{2sl+\nu }xe_i\ot e_j)\\
&=&\sum _{l, t, j=0}^{2\nu -1}
\omega ^{2slt+2stj}e_l\ot e_t\ot e_j+\beta 
(\sum _{l, t=0}^{2\nu -1}
g^{2st}e_l\ot e_tx\ot g^{2st+\nu }x \\
&&-\sum _{t, i=0}^{2\nu -1}
g^{2st}\ot e_txe_i\ot g^{2si+\nu }x
+\sum _{l, i, j=0}^{2\nu -1}
e_lx\ot g^{2sl+2sj+\nu }xe_i\ot e_j\\
&&-\sum _{l, i, j=0}^{2\nu -1}
xe_l\ot g^{2sl+2sj+\nu }xe_i\ot e_j)\\
&=&\sum _{l, t, j=0}^{2\nu -1}
\omega ^{2slt+2stj}e_l\ot e_t\ot e_j+\beta 
(\sum _{t=0}^{2\nu -1}
g^{2st}\ot e_tx\ot g^{2st+\nu }x \\
&&-\sum _{t, i=0}^{2\nu -1}
g^{2st}\ot e_txe_i\ot g^{2si+\nu }x
+\sum _{l, j=0}^{2\nu -1}
e_lx\ot g^{2sl+2sj+\nu }x\ot e_j\\
&&-\sum _{l, j=0}^{2\nu -1}
xe_l\ot g^{2sl+2sj+\nu }x\ot e_j), 
\end{eqnarray*}

\begin{eqnarray*}
T_{23}T_{12}&=&\sum _{l, t, j=0}^{2\nu -1}
\omega ^{2slt+2stj}e_l\ot e_t\ot e_j+\beta 
(\sum _{l, t, j=0}^{2\nu -1}(-1)^t
\omega ^{2slt+2stj}e_lx\ot e_tx\ot e_j \\
&&-\sum _{l, t, j=0}^{2\nu -1}(-1)^{t}
\omega ^{2slt+2stj}xe_l\ot e_tx\ot e_j
+\sum _{l, t, i, j=0}^{2\nu -1}(-1)^j
\omega ^{2slt+2sij}e_l\ot e_ixe_t\ot e_jx\\
&&-\sum _{l, t, j=0}^{2\nu -1}(-1)^{j}
\omega ^{2slt+2stj}e_l\ot xe_t\ot e_jx)\\
&=&\sum _{l, t, j=0}^{2\nu -1}
\omega ^{2slt+2stj}e_l\ot e_t\ot e_j+\beta 
(\sum _{l, t, j=0}^{2\nu -1}(-1)^t
e_lx\ot g^{2sl+2sj}e_tx\ot e_j \\
&&-\sum _{l, t, j=0}^{2\nu -1}(-1)^{t}
xe_l\ot g^{2sl+2sj}e_tx\ot e_j
+\sum _{l, t, i, j=0}^{2\nu -1}(-1)^j
\omega ^{2sli+2stj}e_l\ot e_txe_i\ot e_jx\\
&&-\sum _{l, t, j=0}^{2\nu -1}(-1)^{j}
\omega ^{2slt}e_l\ot xe_t\ot g^{2st}e_jx)\\
&=&\sum _{l, t, j=0}^{2\nu -1}
\omega ^{2slt+2stj}e_l\ot e_t\ot e_j+\beta 
(\sum _{l, j=0}^{2\nu -1}
e_lx\ot g^{2sl+2sj+\nu }x\ot e_j \\
&&-\sum _{l, j=0}^{2\nu -1}
xe_l\ot g^{2sl+2sj+\nu }x\ot e_j
+\sum _{l, t, i, j=0}^{2\nu -1}(-1)^j
g^{2si}e_l\ot e_txe_i\ot g^{2st}e_jx\\
&&-\sum _{l, t=0}^{2\nu -1}
g^{2st}e_l\ot xe_t\ot g^{2st+\nu }x)\\
&=&\sum _{l, t, j=0}^{2\nu -1}
\omega ^{2slt+2stj}e_l\ot e_t\ot e_j+\beta 
(\sum _{l, j=0}^{2\nu -1}
e_lx\ot g^{2sl+2sj+\nu }x\ot e_j \\
&&-\sum _{l, j=0}^{2\nu -1}
xe_l\ot g^{2sl+2sj+\nu }x\ot e_j
+\sum _{t, i=0}^{2\nu -1}
g^{2si}\ot e_txe_i\ot g^{2st+\nu }x\\
&&-\sum _{t=0}^{2\nu -1}
g^{2st}\ot xe_t\ot g^{2st+\nu }x).
\end{eqnarray*}
Thus, we can see that we have 
\begin{eqnarray*}
T_{12}T_{23}-T_{23}T_{12}&=&\beta 
(\sum _{t=0}^{2\nu -1}
g^{2st}\ot e_tx\ot g^{2st+\nu }x 
-\sum _{t, i=0}^{2\nu -1}
g^{2st}\ot e_txe_i\ot g^{2si+\nu }x\\
&&-\sum _{t, i=0}^{2\nu -1}
g^{2si}\ot e_txe_i\ot g^{2st+\nu }x
+\sum _{t=0}^{2\nu -1}
g^{2st}\ot xe_t\ot g^{2st+\nu }x).
\end{eqnarray*}
We need to prove now that we have 
\begin{eqnarray*} 
&&xe_{l}=e_{l-\nu}x, 
\end{eqnarray*}
for all $0\leq l\leq 2\nu -1$, where the subscripts are taken mod $2\nu $.
We use the following facts:
\begin{eqnarray*}
&&\omega^{\nu}=-1,\\
&&xg^i=(-1)^{i}g^ix=\omega^{i\nu}g^ix.
\end{eqnarray*}
We have:
\begin{eqnarray*} 
xe_l&=&x\frac{1}{2\nu }\displaystyle\sum _{i=0}^{2\nu -1}\omega ^{-il}g^i\\
&=&\frac{1}{2\nu }\displaystyle\sum _{i=0}^{2\nu -1}\omega ^{-il}xg^i\\
&=&\frac{1}{2\nu }\displaystyle\sum _{i=0}^{2\nu -1}\omega ^{-il}\omega^{i\nu}g^ix\\
&=&\frac{1}{2\nu }\displaystyle\sum _{i=0}^{2\nu -1}\omega ^{-il+i\nu}g^ix\\
&=&\frac{1}{2\nu }\displaystyle\sum _{i=0}^{2\nu -1}\omega ^{-i(l-\nu)}g^ix\\
&=&e_{l-\nu}x, \;\;\;q.e.d.
\end{eqnarray*}
Now we compute:
\begin{eqnarray*}
\sum _{t, i=0}^{2\nu -1}g^{2st}\ot e_txe_i\ot g^{2si+\nu }x
&=&\sum _{t, i=0}^{2\nu -1}g^{2st}\ot e_te_{i-\nu}x\ot g^{2si+\nu }x\\
&=&\sum _{t, i=0}^{2\nu -1}g^{2st}\ot \delta_{t,i-\nu}e_tx\ot g^{2si+\nu }x\\
&=&\sum _{t=0}^{2\nu -1}g^{2st}\ot e_tx\ot g^{2s(t+\nu)+\nu }x\\
&=&\sum _{t=0}^{2\nu -1}g^{2st}\ot e_tx\ot g^{2st+\nu}g^{2s\nu}x\\
&=&\sum _{t=0}^{2\nu -1}g^{2st}\ot e_tx\ot g^{2st+\nu}x, 
\end{eqnarray*}
so we have $\displaystyle\sum _{t=0}^{2\nu -1}g^{2st}\ot e_tx\ot g^{2st+\nu}x-
\displaystyle\sum _{t, i=0}^{2\nu -1}g^{2st}\ot e_txe_i\ot g^{2si+\nu }x=0$. 
Similarly, we have:
\begin{eqnarray*}
\sum _{t, i=0}^{2\nu -1}g^{2si}\ot e_txe_i\ot g^{2st+\nu }x
&=&\sum _{t, i=0}^{2\nu -1}g^{2si}\ot xe_{t+\nu}e_i\ot g^{2st+\nu }x\\
&=&\sum _{t, i=0}^{2\nu -1}g^{2si}\ot x\delta_{t+\nu,i}e_i\ot g^{2st+\nu }x\\
&=&\sum _{i=0}^{2\nu -1}g^{2si}\ot xe_i\ot g^{2s(i-\nu)+\nu }x\\
&=&\sum _{i=0}^{2\nu -1}g^{2si}\ot xe_i\ot g^{2si+\nu }x, 
\end{eqnarray*}
so we have $\displaystyle\sum _{t=0}^{2\nu -1}g^{2st}\ot xe_t\ot g^{2st+\nu }x-
\displaystyle\sum _{t, i=0}^{2\nu -1}g^{2si}\ot e_txe_i\ot g^{2st+\nu }x=0$.  
Consequently, we have $T_{12}T_{23}-T_{23}T_{12}=0$, and so we 
obtained:
\begin{theorem}
Any quasitriangular structure $R_{s, \beta }$ on Radford's Hopf algebra $H_{\nu }$ is pseudotriangular. 
\end{theorem}
%%%%%%%%%%%%%%%%%%%%%%%%%%%%%%%%
\section{Hopf algebras with positive bases}
%%%%%%%%%%%%%%%%%%%%%%%%%%%%%%%%%%%
${\;\;\;\;\;}$In this section the base field is assumed to be 
${\mathbb C}$, the field of complex numbers. 

We recall from \cite{LYZ1} that a basis of a Hopf algebra over 
${\mathbb C}$ is called {\em positive} if all the structure constants 
(for the unit, counit, multiplication, comultiplication and antipode) 
with respect to this basis are nonnegative real numbers. Also, a 
quasitriangular structure $R$ on a Hopf algebra having a positive basis $B$ 
is called {\em positive} in \cite{LYZ2} if the coefficients of $R$ in the 
basis $B\otimes B$ are nonnegative real numbers. The finite 
dimensional Hopf algebras having a positive basis and the positive quasitriangular 
structures on them have been classified in \cite{LYZ1}, 
\cite{LYZ2} as follows. 

Let $G$ be a group (we denote by $e$ its unit). A {\em unique 
factorization} $G=G_+G_-$ of $G$ consists of two subgroups 
$G_+$ and $G_-$ of $G$ such that any $g\in G$ can be written uniquely 
as $g=g_+g_-$, with $g_+\in G_+$ and $g_-\in G_-$. By considering the inverse map, 
we can also write uniquely $g=\overline{g}_-
\overline{g}_+$, with $\overline{g}_-\in G_-$ and 
$\overline{g}_+\in G_+$. 

Let $u\in G_+$, $x\in G_-$; then we can write uniquely
\begin{eqnarray*}
&&xu=(^xu)(x^u), \;\;with\;\;^xu\in G_+\;\;and \;\;x^u\in G_-, \\
&&ux=(^ux)(u^x), \;\;with\;\;^ux\in G_- \;\;and\;\;u^x\in G_+. 
\end{eqnarray*}
So, we have the following actions of $G_+$ and $G_-$ on each other 
(from left and right):
\begin{eqnarray*}
&&G_-\times G_+\rightarrow G_+, \;\;\;(x, u)\mapsto \;^xu, \\
&&G_-\times G_+\rightarrow G_-, \;\;\;(x, u)\mapsto \;x^u, \\
&&G_+\times G_-\rightarrow G_-, \;\;\;(u, x)\mapsto \;^ux, \\
&&G_+\times G_-\rightarrow G_+, \;\;\;(u, x)\mapsto \;u^x.
\end{eqnarray*}
The relations between these actions and the decompositions 
$g=g_+g_-=\overline{g}_-\overline{g}_+$ are: \\
$^{\overline{g}_-}\overline{g}_+=g_+$; $\overline{g}_-^
{\overline{g}_+}=g_-$; $g_+^{g_-}=\overline{g}_+$; 
$^{g_+}g_-=\overline{g}_-$; $(^{g_+}g_-)(g_+^{g_-})=g_+g_-$; 
$(^{g_-}g_+)(g_-^{g_+})=g_-g_+$. 

Given a unique factorization $G=G_+G_-$ of a finite group $G$, one can 
construct a finite dimensional Hopf algebra $H(G; G_+, G_-)$, which is 
the vector space spanned by the set $G$ (we denote by $\{g\}$ an 
element $g\in G$ when it is regarded as an element in 
$H(G; G_+, G_-)$) with the following Hopf algebra structure:\\
multiplication: $\;\;\;\{g\}\{h\}=\delta _{g_+^{g_-}, h_+}\{gh_-\}$\\
unit: $\;\;\;1=\sum _{g_+\in G_+}\{g_+\}$\\
comultiplication: $\;\;\;\Delta (\{g\})=\sum _{h_+\in G_+}
\{g_+h_+^{-1}(^{h_+}g_-)\}\otimes \{h_+g_-\}$\\
counit: $\;\;\;\varepsilon (\{g\})=\delta _{g_+, e}$\\
antipode: $\;\;\;S(\{g\})=\{g^{-1}\}$

The Hopf algebra $H(G; G_+, G_-)$ has $G$ as the obvious positive basis. 
Conversely, it was proved in \cite{LYZ1} that all finite dimensional Hopf 
algebras with positive bases are of the form $H(G; G_+, G_-)$. 

The positive quasitriangular and triangular structures on $H(G; G_+, G_-)$ 
have been described in \cite{LYZ2} as follows:  
\begin{theorem} (\cite{LYZ2})
Let $G=G_+G_-$ be a unique factorization of a finite group $G$. 
Let $\xi , \eta :G_+\rightarrow G_-$ be two group homomorphisms 
satisfying the following conditions: 
\begin{eqnarray}
&&\xi (u)^v=\xi (u^{\eta (v)}), \label{rel6} \\
&&^u\eta (v)=\eta (^{\xi (u)}v), \label{rel7} \\
&&uv=(^{\xi (u)}v)(u^{\eta (v)}), \label{rel8} \\
&&\xi (^xu)x^u=x\xi (u), \label{rel9} \\
&&\eta (^xu)x^u=x\eta (u), \label{rel10}
\end{eqnarray}
for all $u, v\in G_+$ and $x\in G_-$. Then 
\begin{eqnarray*}
&&R(\xi , \eta ):=\sum _{u, v\in G_+}\{u(\eta (v)^u)^{-1}\}\ot 
\{v\xi (u)\}
\end{eqnarray*}
is a positive quasitriangular structure on $H(G; G_+, G_-)$. 
Conversely, every positive quasitriangular structure on 
$H(G; G_+, G_-)$ is given by the above construction. 

Moreover, each of the conditions (\ref{rel6})-(\ref{rel10}) is 
equivalent to the corresponding property below:
\begin{eqnarray}
&&^v\xi (u)=\xi (^{\eta (v)}u), \label{rel11} \\
&&\eta (v)^u=\eta (v^{\xi (u)}), \label{rel12} \\
&&uv=(^{\eta (u)}v)(u^{\xi (v)}), \label{rel13} \\
&&^ux\xi (u^x)=\xi (u)x, \label{rel14} \\
&&^ux\eta (u^x)=\eta (u)x. \label{rel15}
\end{eqnarray}

Moreover, $R(\xi , \eta )$ is triangular if and only if $\xi =\eta $.
\end{theorem}

Our aim now is to characterize those $R(\xi , \eta )$ that are 
pseudotriangular. So, let $R=R(\xi , \eta )$ be a positive quasitriangular 
structure on 
$H(G; G_+, G_-)$. We have (see \cite{LYZ2}): 
\begin{eqnarray*}
&&R_{21}R=\sum _{u, v\in G_+}\{v\xi (u)(\eta (\overline{v})
^{\overline{u}})^{-1}\}\ot \{u(\eta (v)^u)^{-1}
\xi (\overline{u})\}, 
\end{eqnarray*}
where we denoted $\overline{u}=v^{\xi (u)}$ and 
$\overline{v}=\;^{\eta (v)}u$. 

We denote $T=R_{21}R$ and we compute (by using the formula 
for the multiplication of $H(G; G_+, G_-))$:
\begin{eqnarray*}
T_{12}T_{23}&=&
(\sum _{u, v\in G_+}\{v\xi (u)(\eta (\overline{v})
^{\overline{u}})^{-1}\}\ot \{u(\eta (v)^u)^{-1}
\xi (\overline{u})\}\ot 1)\\
&&
(\sum _{s, t\in G_+}1\ot \{t\xi (s)(\eta (\overline{t})
^{\overline{s}})^{-1}\}\ot \{s(\eta (t)^s)^{-1}
\xi (\overline{s})\})\\
&=&\sum _{u, v, s, t\in G_+}\{v\xi (u)(\eta (\overline{v})
^{\overline{u}})^{-1}\}\ot \{u(\eta (v)^u)^{-1}
\xi (\overline{u})\}\{t\xi (s)(\eta (\overline{t})
^{\overline{s}})^{-1}\}\\
&&\ot \{s(\eta (t)^s)^{-1}
\xi (\overline{s})\})\\
&=&\sum _{u, v, s\in G_+}\{v\xi (u)(\eta (^{\eta (v)}u)^
{(v^{\xi (u)})})^{-1}\}\ot \{u(\eta (v)^u)^{-1}\xi (v^{\xi (u)})
\xi (s)(\eta (^{\eta (t)}s)^{(t^{\xi (s)})})^{-1}\}\\
&&\ot \{s(\eta (t)^s)^{-1}\xi (t^{\xi (s)})\},
\end{eqnarray*}
where $t=u^{(\eta (v)^u)^{-1}\xi (v^{\xi (u)})}$, and 
\begin{eqnarray*}
T_{23}T_{12}&=&
(\sum _{a, b\in G_+}1\ot \{b\xi (a)(\eta (^{\eta (b)}a)^
{(b^{\xi (a)})})^{-1}\}\ot \{a(\eta (b)^a)^{-1}\xi (b^{\xi (a)})\})\\
&&(\sum _{c, d\in G_+}\{d\xi (c)(\eta (^{\eta (d)}c)^
{(d^{\xi (c)})})^{-1}\}\ot \{c(\eta (d)^c)^{-1}\xi (d^{\xi (c)})\}\ot 1)\\
&=&\sum _{a, b, c, d\in G_+}
\{d\xi (c)(\eta (^{\eta (d)}c)^
{(d^{\xi (c)})})^{-1}\}\ot \{b\xi (a)(\eta (^{\eta (b)}a)^
{(b^{\xi (a)})})^{-1}\}
\{c(\eta (d)^c)^{-1}\xi (d^{\xi (c)})\}\\
&&\ot 
\{a(\eta (b)^a)^{-1}\xi (b^{\xi (a)})\}\\
&=&\sum _{a, b, d\in G_+}
\{d\xi (c)(\eta (^{\eta (d)}c)^
{(d^{\xi (c)})})^{-1}\}
\ot \{b\xi (a)(\eta (^{\eta (b)}a)^
{(b^{\xi (a)})})^{-1}
(\eta (d)^c)^{-1}\xi (d^{\xi (c)})\}\\
&&\ot 
\{a(\eta (b)^a)^{-1}\xi (b^{\xi (a)})\},
\end{eqnarray*}
where $c=b^{\xi (a)(\eta (^{\eta (b)}a)^
{(b^{\xi (a)})})^{-1}}$. By writing down what means $T_{12}T_{23}
=T_{23}T_{12}$, we obtain:
\begin{proposition}
The positive quasitriangular structure $R(\xi , \eta )$ is 
pseudotriangular if and only if the following conditions are satisfied:
\begin{eqnarray*}
&&\xi (u)(\eta (^{\eta (v)}u)^
{(v^{\xi (u)})})^{-1}=\xi (c)(\eta (^{\eta (v)}c)^
{(v^{\xi (c)})})^{-1}, \\
&&(\eta (v)^u)^{-1}\xi (v^{\xi (u)})
\xi (s)(\eta (^{\eta (t)}s)^{(t^{\xi (s)})})^{-1}=
\xi (s)(\eta (^{\eta (u)}s)^
{(u^{\xi (s)})})^{-1}
(\eta (v)^c)^{-1}\xi (v^{\xi (c)}), \\
&&(\eta (t)^s)^{-1}\xi (t^{\xi (s)})=
(\eta (u)^s)^{-1}\xi (u^{\xi (s)}), 
\end{eqnarray*}
for all $u, v, s\in G_+$, where 
$t=u^{(\eta (v)^u)^{-1}\xi (v^{\xi (u)})}$ and 
$c=u^{\xi (s)(\eta (^{\eta (u)}s)^
{(u^{\xi (s)})})^{-1}}$.
\end{proposition}

A better description may be obtained for a certain class of positive 
quasitriangular structures. 
\begin{definition}(\cite{LYZ2})
A positive quasitriangular structure $R(\xi , \eta )$ on 
$H(G; G_+, G_-)$ is called {\em normal} if $\xi (u)=e$ for all 
$u\in G_+$. 
\end{definition}
\begin{theorem}\label{normal}
A normal positive quasitriangular structure $R(\xi , \eta )$ on 
$H(G; G_+, G_-)$ is pseudotriangular if and only if $\eta (uv)=
\eta (vu)$ for all $u, v\in G_+$.  
\end{theorem}
\begin{proof}
We note first that, since $\xi (u)=e$ for all $u\in G_+$, some of the 
relations (\ref{rel6})-(\ref{rel15}) may be simplified, in particular we have 
$^u\eta (v)=\eta (v)$, 
$uv=v(u^{\eta (v)})$, 
$\eta (v)^u=\eta (v)$, 
$uv=(^{\eta (u)}v)u$, 
for all $u, v\in G_+$. By using these relations, 
together with the fact that $\xi (u)=e$ for all $u\in G_+$, the 
three conditions in the above Proposition may be also simplified, so 
we obtain that $R(\xi , \eta )$ is pseudotriangular if and only if we have:
\begin{eqnarray*}
&&\eta (vuv^{-1})=\eta (vcv^{-1}), \\
&&\eta (v)^{-1}\eta (tst^{-1})^{-1}=
\eta (usu^{-1})^{-1}\eta (v)^{-1}, \\
&&\eta (t)^{-1}=\eta (u)^{-1},
\end{eqnarray*} 
for all $u, v, s\in G_+$, where $t=vuv^{-1}$ and 
$c=usus^{-1}u^{-1}$, and one can easily see that each of these 
three conditions is equivalent to the condition $\eta (uv)=\eta (vu)$, 
for all $u, v\in G_+$. 
\end{proof}

We recall from \cite{psvo} that the canonical quasitriangular structure on the 
Drinfeld double of a finite dimensional Hopf algebra $H$ is 
pseudotriangular if and only if $H$ is commutative and cocommutative. 
In particular, if $G$ is a finite group, the canonical quasitriangular structure 
on the Drinfeld double of the dual $k[G]^*$ of the group algebra 
$k[G]$ is pseudotriangular if and only if $G$ is abelian. We 
want to reobtain this result (over ${\mathbb C}$) as an application 
of Theorem \ref{normal}.

We consider the unique factorization $G=G_+G_-$, where 
$G_+=G$ and $G_-=\{e\}$ (so the Hopf algebra $H(G; G_+, G_-)$ 
is exactly $k[G]^*$). As in \cite{LYZ2}, we consider the group 
$\tilde{G}=G\times G$, with the unique factorization $\tilde{G}=
\tilde{G}_+\tilde{G}_-$, where $\tilde{G}_+=G\times \{e\}$ and 
$\tilde{G}_-=\{(g, g): g\in G\}$.  By \cite{LYZ2}, the group 
homomorphisms $\xi , \eta :\tilde{G}_+\rightarrow \tilde{G}_-$ 
defined by $\xi (g, e)=(e, e)$ and $\eta (g, e)=(g, g)$ induce a 
positive quasitriangular structure $R(\xi , \eta )$ on 
$H(\tilde{G}; \tilde{G}_+, \tilde{G}_-)$ and moreover 
$H(\tilde{G}; \tilde{G}_+, \tilde{G}_-)$ is the Drinfeld double 
of $H(G; G_+, G_-)=k[G]^*$ and $R(\xi , \eta )$ is its canonical 
quasitriangular structure. Obviously $R(\xi , \eta )$ is normal, so  
we may apply Theorem \ref{normal} and we obtain that $R(\xi , \eta )$ 
is pseudotriangular if and only if $(gh, gh)=(hg, hg)$ for all 
$g, h\in G$, i.e. if and only if $G$ is abelian.
%%%%%%%%%%%%%%%%%%%%%%%%%%%%%%%%%
\section{Universality of the pseudosymmetric category ${\cal PS}$ }
%%%%%%%%%%%%%%%%%%%%%%%%%%%%%%%%
\setcounter{equation}{0}
%%%%%%%%%%%%%%%%%%%%%%%%%%%%%%
${\;\;\;\;}$
In this section we use terminology, notation and some results from \cite{k} (but we 
use the term ''monoidal'' instead of ''tensor'' when we speak about tensor categories 
and tensor functors).  

Our aim is to show that the pseudosymmetric category ${\cal PS}$ introduced in 
\cite{ps} has two universality properties similar to the ones of the 
braid category ${\cal B}$, the universal braided monoidal category (see \cite{k}). 
First, we recall  from \cite{ps} the definition of ${\cal PS}$. 
The objects of ${\cal PS}$ are natural numbers $n\in  {\mathbb N}$. 
The set of morphisms from 
$m$ to $n$ is empty if $m\neq n$ and is $PS_n:=\frac{B_n}{[P_n,P_n]}$ if $m = n$, 
where $B_n$ (respectively $P_n$) is the braid group (respectively pure braid 
group) on $n$ strands. 
The monoidal structure of ${\cal PS}$ is defined as the one for ${\cal B}$, 
and so is the braiding, namely (we denote as usual by $\sigma _1, \sigma _2, ... , 
\sigma _{n-1}$ the standard generators of $B_n$ and by $\pi_n$ the natural morphism from  $B_n$ to $PS_n$):
\begin{eqnarray*}
&c_{n,m} : n \otimes  m \to m \otimes n, \;\;\;c_{0,n} = id_n = c_{n,0},&\\
&c_{n,m} = \pi_{n+m}((\sigma_{m}\sigma_{m-1}\cdots \sigma_1) 
(\sigma_{m+1}\sigma_{m}\cdots \sigma_2)\cdots 
(\sigma_{m+n-1}\sigma_{m+n-2}\cdots \sigma_{n})) \;\; {\rm if}\;  m, n > 0.&
\end{eqnarray*}

In order to introduce the first universality property for ${\cal PS}$, we need the following  definition, motivated by results in \cite{ps} and by the definition of Yang-Baxter operators from \cite{k}:
\begin{definition} If $V$ is an object in a monoidal category 
$({\cal C}, \otimes, I, a, l,r)$, an automorphism $\sigma$ of $V\otimes V$ is called a {\em pseudosymmetric Yang-Baxter operator} on $V$ if the following two dodecagons 
(for $\sigma$ and $\sigma^{-1}$) commute:
\begin{displaymath}
    \xymatrix{&\ar[dl]_{\sigma\otimes id_V} \; \; \; \; \; \; \; \; \; \; (V\otimes V)\otimes V 
\; \; \; \; \; \; \; \; \; \;   \ar[dr]^{a_{V,V,V}}&\\
       (V\otimes V)\otimes V \ar[d]_{a_{V,V,V}}& & V\otimes (V\otimes V)\ar[d]^{id_V\otimes \sigma}\\
 V\otimes (V\otimes V) \ar[d]_{id_V\otimes \sigma^{\pm 1}}& & V\otimes (V\otimes V)\ar[d]^{a_{V,V,V}^{-1}}\\
 V\otimes (V\otimes V) \ar[d]_{a_{V,V,V}^{-1}}& & (V\otimes V)\otimes V\ar[d]^{\sigma^{\pm 1}\otimes id_V}\\
 (V\otimes V)\otimes V \ar[d]_{\sigma\otimes id_V}& & (V\otimes V)\otimes V\ar[d]^{a_{V,V,V}}\\
 (V\otimes V)\otimes V \ar[dr]_{a_{V,V,V}}& & V\otimes (V\otimes V)\ar[dl]^{id_V\otimes \sigma}\\
 &V\otimes (V\otimes V)&                    }
\end{displaymath}
\end{definition}

Note that a pseudosymmetric Yang-Baxter operator is a special type of  Yang-Baxter 
operator as defined in  \cite{k}, p. 323. Moreover, just like Yang-Baxter operators, 
they can be transferred by using functors between monoidal categories: 
\begin{lemma} Let $(F,\varphi_0,\varphi_2):{\cal C}\to {\cal D}$ be a monoidal 
functor between two monoidal categories. 
If $\sigma\in Aut(V\otimes V)$ is a pseudosymmetric Yang-Baxter operator on the 
object $V\in {\cal C}$, then 
\begin{eqnarray*}
\sigma'=\varphi_2(V,V)^{-1}\circ F(\sigma)\circ \varphi_2(V,V)
\end{eqnarray*}
is a pseudosymmetric Yang-Baxter operator on $F(V)$. \label{lemma1}
\end{lemma}
\begin{proof}  The proof follows exactly  as in  \cite{k}, Lemma XIII.3.2, 
by using also the identity
\begin{eqnarray*}
&&(\sigma')^{-1}= \varphi_2(V,V)^{-1}\circ F(\sigma^{-1})\circ \varphi_2(V,V)
\end{eqnarray*}
in order to prove the pseudosymmetry of $\sigma '$. 
\end{proof}

We define the  category  $PSYB({\cal C})$ of pseudosymmetric Yang-Baxter operators 
to be a full subcategory of $YB({\cal C})$, the category of Yang-Baxter operators 
defined in \cite{k}.  
An object in $PSYB({\cal C})$ is  a pair $(V,\sigma)$ where $V$ is a object in ${\cal C}$ and $\sigma$ is a pseudosymmetric Yang-Baxter operator. 
%Notice that $(V,id_{V\otimes V})$ is in $PSYB({\cal C})$.

Recall the following construction from \cite{k}. 
Suppose that $(F,\varphi_0,\varphi_2):{\cal B}\to {\cal C}$ is a monoidal functor from the 
universal braid category ${\cal B}$ to a given monoidal category ${\cal C}$. 
Since $c_{1,1}=\sigma_1$ is a Yang-Baxter operator on the object $1\in {\cal B}$, 
it follows that 
$\sigma=\varphi_2^{-1}(1,1)F(c_{1,1})\varphi_2(1,1)$
is a Yang-Baxter operator on $F(1)\in {\cal C}$. In this way we get a functor 
$\Theta:Tens({\cal B},{\cal C})\to YB({\cal C})$, where $Tens({\cal B},{\cal C})$ is the 
category of monoidal functors from ${\cal B}$ to ${\cal C}$. 
It was proved in  \cite{k} that:
\begin{theorem} (\cite{k})
For any monoidal category ${\cal C}$, the functor 
$\Theta:Tens({\cal B},{\cal C})\to YB({\cal C})$ is an equivalence of categories. 
\end{theorem}

One can note that we have a natural monoidal functor $\pi: {\cal B}\to {\cal PS}$ 
induced by the group epimorphism 
$\pi_n:B_n\to PS_n$. This allows us to identify the category $Tens({\cal PS},{\cal C})$ 
with a subcategory of $Tens({\cal B},{\cal C})$. More precisely, we identify it with the full subcategory of all monoidal functors  $F:{\cal B}\to {\cal C}$ with the property that there 
exists a monoidal functor $G:{\cal PS}\to {\cal C}$ such that $F=G\circ\pi$. 

We can state now the first universality property of  ${\cal PS}$:
\begin{theorem}
For any monoidal category ${\cal C}$, the functor 
$\widetilde{\Theta}: Tens({\cal PS},{\cal C})\to PSYB({\cal C})$, $\widetilde{\Theta}(G)=
\Theta(G\circ\pi)$ is an equivalence of categories.
\end{theorem}
\begin{proof} First we note that $\pi(c_{1,1})$ is a pseudosymmetric 
Yang-Baxter operator in ${\cal PS}$ and so by Lemma \ref{lemma1} we have 
$\varphi_2^{-1}(1,1)G(\pi(c_{1,1}))\varphi_2(1,1)\in  PSYB({\cal C})$. This means that $\widetilde{\Theta}$ is well defined. 
Since $\Theta$ is fully faithful  and $\widetilde{\Theta}$ is its restriction to a full subcategory, 
it is enough to show that $\widetilde{\Theta}$ is essentially surjective. 
This follows from the next lemma.
\end{proof}
\begin{lemma} Let ${\cal C}$ be a strict monoidal category and $(V,\sigma)$ an object in $PSYB({\cal C})$. Then there exists a unique strict monoidal functor 
$G:{\cal PS}\to {\cal C}$ such that $G(1)=V$ and $G(\pi(c_{1,1}))=\sigma$. 
\end{lemma}
\begin{proof} From \cite{k}, Lemma XIII.3.5 we know that for all 
$(V,\sigma)\in YB({\cal C})$ there exists a unique strict monoidal functor 
$F:{\cal B}\to {\cal C}$ such that $F(1)=V$ and $F(c_{1,1})=\sigma$. 
It is enough to show that when $(V,\sigma)\in PSYB({\cal C})$ the functor $F$ factors 
through $\pi$. But this follows immediately from the fact (see \cite{ps}) that 
\begin{eqnarray*}
&&PS_n=\frac{B_n}{<\sigma_i\sigma_{i+1}^{-1}\sigma_i=\sigma_{i+1}\sigma_{i}^{-1}\sigma_{i+1}\vert 
1\leq i\leq n-2>}
\end{eqnarray*}
 and the definition of a pseudosymmetric Yang-Baxter operator.
\end{proof}

\begin{definition}(\cite{k}) A monoidal functor $(F,\varphi_0,\varphi _2)$ from a 
braided monoidal category ${\cal C}$ to a braided monoidal category ${\cal D}$ 
is braided if for every pair $(U,V)$ of objects in ${\cal C}$ the square
%$$F(c_{U,V})\varphi_2(U,V)=\varphi_2(V,U)c_{F(U),F(V)}$$
\begin{displaymath}
    \xymatrix{    &F(U)\otimes F(V)\ar[d]_{c_{F(U),F(V)}} \ar[r]^{\varphi_2} &F(U\otimes V)\ar[d]^{F(c_{U,V})}\\
                   &F(V)\otimes F(U) \ar[r]^{\varphi_2} &F(V\otimes U)                  }
\end{displaymath}
commutes. Denote by $Br({\cal C},{\cal D})$ the category  whose objects are braided 
monoidal functors and morphisms are natural monoidal transformations. 
\end{definition}
\begin{theorem} (\cite{k}) \label{5.7}
For a braided monoidal category ${\cal C}$, the functor 
$\Theta':Br({\cal B},{\cal C})\to {\cal C}$ defined by $\Theta'(F)=F(1)$ is an 
equivalence of categories. 
\end{theorem}

In the definition of a pseudosymmetric braided category ${\cal C}$ 
introduced in \cite{psvo} was assumed that ${\cal C}$ was a {\em strict} monoidal category.   
The next proposition is the analogue of Theorem 3.7  from \cite{psvo} for monoidal 
categories with nontrivial associativity constraints.  
Note that the proof that we present here is very direct and is inspired by the 
results in \cite{ps}.  
\begin{proposition} Let $({\cal C},\otimes,I,a,l,r,c)$ be a braided monoidal category. 
The following conditions are equivalent:\\
(i) For every $U$,  $V$, $W\in {\cal C}$ the following diagram is commutative: 
\begin{displaymath}
    \xymatrix{&\ar[dl]_{c_{U,V}\otimes id_W} \; \; \; \; \; \; \; \; \; \; (U\otimes V)\otimes W 
\; \; \; \; \; \; \; \; \; \;   \ar[dr]^{a_{U,V,W}}&\\
       (V\otimes U)\otimes W \ar[d]_{a_{V,U,W}}& & U\otimes (V\otimes W)\ar[d]^{id_U\otimes c_{V,W}}\\
 V\otimes (U\otimes W) \ar[d]_{id_V\otimes c_{W,U}^{-1}}& & U\otimes (W\otimes V)\ar[d]^{a_{U,W,V}^{-1}}\\
 V\otimes (W\otimes U) \ar[d]_{a_{V,W,U}^{-1}}& & (U\otimes W)\otimes V\ar[d]^{c_{W,U}^{-1}\otimes id_V}\\
 (V\otimes W)\otimes U \ar[d]_{c_{V,W}\otimes id_U}& & (W\otimes U)\otimes V\ar[d]^{a_{W,U,V}}\\
 (W\otimes V)\otimes U \ar[dr]_{a_{W,V,U}}& & W\otimes (U\otimes V)\ar[dl]^{id_W\otimes c_{U,V}}\\
 &W\otimes (V\otimes U)&                  }\label{diag1}  
\end{displaymath}
(ii) For every $U$,  $V$, $W\in {\cal C}$ the following diagram is commutative: 
\begin{displaymath}
    \xymatrix{&\ar[dl]_{c_{V,U}c_{U,V}\otimes id_W} \; \; \; \; \; \; \; \; \; \; 
(U\otimes V)\otimes W \; \; \; \; \; \; \; \; \; \;   \ar[dr]^{a_{U,V,W}}&\\
       (U\otimes V)\otimes W \ar[d]_{a_{U,V,W}}& & U\otimes (V\otimes W)\ar[d]^{id_U\otimes c_{W,V}c_{V,W}}\\
 U\otimes (V\otimes W) \ar[d]_{id_U\otimes c_{W,V}c_{V,W}}& & U\otimes (V\otimes W)\ar[d]^{a_{U,V,W}^{-1}}\\
 U\otimes (V\otimes W) \ar[dr]_{a_{U,V,W}^{-1}}& & (U\otimes V)\otimes W\ar[dl]^{c_{V,U}c_{U,V}\otimes id_W}\\
 &(U\otimes V)\otimes W&                    }
\end{displaymath} 
\label{prop2}
\end{proposition}
\begin{proof} Take  $U$,  $V$, $W\in {\cal C}$. Using only the fact that ${\cal C}$ is a braided category we have\\[2mm]
${\;\;\;\;\;}$$((c_{V,U}c_{U,V})\otimes id_W)a_{U,V,W}^{-1}(id_U\otimes (c_{W,V} c_{V,W}))a_{U,V,W}$
\begin{eqnarray*}
&=&(c_{V,U}\otimes id_W)a_{V,U,W}^{-1}(id_V\otimes c_{U,W}^{-1})[(id_V\otimes c_{U,W})a_{V,U,W}(c_{U,V}\otimes id_W)]\\&&a_{U,V,W}^{-1}(id_U\otimes c_{W,V} c_{V,W})a_{U,V,W}\\
&=&(c_{V,U}\otimes id_W)a_{V,U,W}^{-1}(id_V\otimes c_{U,W}^{-1})
[a_{V,W,U}c_{U,V\otimes W}a_{U,V,W}]a_{U,V,W}^{-1}\\&&(id_U\otimes c_{W,V}\circ c_{V,W})a_{U,V,W}\\
&=&(c_{V,U}\otimes id_W)a_{V,U,W}^{-1}(id_V\otimes 
c_{U,W}^{-1})a_{V,W,U} (c_{W,V}\otimes id_U)(c_{V,W}\otimes id_U) 
c_{U,V\otimes W} a_{U,V,W}, 
\end{eqnarray*}
${\;\;\;}$$
a_{U,V,W}^{-1}(id_U\otimes c_{W,V})(id_U\otimes c_{V,W})a_{U,V,W}
(c_{V,U}\otimes id_{W})(c_{U,V}\otimes id_W)$
\begin{eqnarray*}
&=&a_{U,V,W}^{-1}(id_U\otimes c_{W,V})a_{U,W,V}c_{V,U\otimes W}a_{V,U,W}
(c_{U,V}\otimes id_W)\\
&=&a_{U,V,W}^{-1}(id_U\otimes c_{W,V})a_{U,W,V}c_{V,U\otimes W}
(id_V\otimes c_{U,W}^{-1})(id_V\otimes c_{U,W})a_{V,U,W}(c_{U,V}\otimes id_W)\\
&=&a_{U,V,W}^{-1}(id_U\otimes c_{W,V})a_{U,W,V}(c_{U,W}^{-1}\otimes id_V)
c_{V,W\otimes U}(id_V\otimes c_{U,W})a_{V,U,W}(c_{U,V}\otimes id_W)\\
&=&a_{U,V,W}^{-1}(id_U\otimes c_{W,V})a_{U,W,V}(c_{U,W}^{-1}\otimes id_V)
a_{W,U,V}^{-1}(id_W\otimes c_{V,U})a_{W,V,U}(c_{V,W}\otimes id_U)a_{V,W,U}^{-1}\\
&&(id_V\otimes c_{U,W})a_{V,U,W}(c_{U,V}\otimes id_W)\\
&=&a_{U,V,W}^{-1}(id_U\otimes c_{W,V})a_{U,W,V}(c_{U,W}^{-1}\otimes id_V)
a_{W,U,V}^{-1}(id_W\otimes c_{V,U})a_{W,V,U}(c_{V,W}\otimes id_U)a_{V,W,U}^{-1}\\
&&a_{V,W,U}c_{U,V\otimes W}a_{U,V,W}\\
&=&a_{U,V,W}^{-1}(id_U\otimes c_{W,V})a_{U,W,V}(c_{U,W}^{-1}\otimes id_V)
a_{W,U,V}^{-1}(id_W\otimes c_{V,U})a_{W,V,U}\\
&&(c_{V,W}\otimes id_U)c_{U,V\otimes W}a_{U,V,W}.
\end{eqnarray*}
This means that the condition (ii) holds if and only if 
\begin{eqnarray*}
&&(c_{V,U}\otimes id_W)a_{V,U,W}^{-1}(id_V\otimes c_{U,W}^{-1})a_{V,W,U} 
(c_{W,V}\otimes id_U) \\
&&\;\;\;\;\;\;\;\;\;\;\;\;\;=
a_{U,V,W}^{-1}(id_U\otimes c_{W,V})a_{U,W,V}
(c_{U,W}^{-1}\otimes id_V)a_{W,U,V}^{-1}(id_W\otimes c_{V,U})a_{W,V,U},
\end{eqnarray*}
and this condition is obviously equivalent with (i).
\end{proof}
\begin{definition} We say that a braided monoidal category $({\cal C},\otimes,I,a,l,r,c)$ is pseudosymmetric if it satisfies any of the two 
equivalent conditions from Proposition \ref{prop2}.
\end{definition}
\begin{remark}
If $({\cal C},\otimes,I,a,l,r,c)$ is a pseudosymmetric braided monoidal category and $V$ is an 
object in ${\cal C}$, then $c_{V, V}$ is a pseudosymmetric Yang-Baxter operator on $V$.
\end{remark}
\begin{lemma}
If the braided category ${\cal C}$ is pseudosymmetric then 
$Br({\cal B},{\cal C})\cong Br({\cal PS},{\cal C})$. 
\end{lemma}
\begin{proof}
The isomorphism  is induced by $\pi:{\cal B}\to {\cal PS}$. 
More precisely, we have 
\begin{eqnarray*}
&&\pi^*:Br({\cal PS},{\cal C})\to  Br({\cal B},{\cal C}), \;\;\; \pi^*(G)=G\circ \pi.
\end{eqnarray*} 
Because $\pi_n:B_n\to PS_n$ is surjective and the category ${\cal C}$ is 
pseudosymmetric, any functor $F\in Br({\cal B},{\cal C})$ is of the form 
$F=G\circ \pi$ for some unique $G\in Br({\cal PS},{\cal C})$. 
\end{proof}

As a consequence of this and  Theorem \ref{5.7} we obtain the second 
universality property of ${\cal PS}$: 
\begin{theorem}  For a pseudosymmetric braided category ${\cal C}$, 
the functor $\widetilde{\Theta'}:Br({\cal PS},{\cal C})\to {\cal C}$ defined by 
$\widetilde{\Theta'}(G)=G(1)$ is an equivalence of categories. 
\end{theorem}

%%%%%%%%%%%%%%%%%%%%%%%%%%%%%%%%%%%%


\begin{thebibliography}{99}
%%%%%%%%%%%%%%%%%%%%%%%%%%%%%%%%%%   
\bibitem{brug}
A. Brugui\`{e}res, Double braidings, twists and tangle invariants,   
{\sl J. Pure Appl. Algebra} {\bf 204} (2006), 170--194.

\bibitem{bcp}
D. Bulacu, S. Caenepeel, F. Panaite, Yetter-Drinfeld categories for 
quasi-Hopf algebras,   
{\sl Comm. Algebra} {\bf 34} (2006), 1--35.

\bibitem{cc}
G. Carnovale, J. Cuadra, The Brauer group of some quasitriangular Hopf algebras, 
{\sl J. Algebra} {\bf 259} (2003), 512--532. 

\bibitem{js}
A. Joyal, R. Street, Braided tensor categories, {\sl Adv. Math.} 
{\bf 102} (1993), 20--78.   

\bibitem{k}
C. Kassel, {\sl Quantum groups}, Grad. Texts Math.  
{\bf 155}, Springer-Verlag, Berlin, 1995.

\bibitem{LYZ1}
J.-H. Lu, M. Yan, Y.-C. Zhu, On Hopf algebras with positive bases, 
{\sl J. Algebra} {\bf 237} (2001), 421--445. 

\bibitem{LYZ2}
J.-H. Lu, M. Yan, Y.-C. Zhu, On quasitriangular structures on Hopf algebras 
with positive bases, in ''New trends in Hopf algebra theory'' (1999), Amer. Math. Soc., 
Providence, RI, pp. 339--356. 

\bibitem{psv}
F. Panaite, M. D. Staic, F. Van Oystaeyen, On some classes of lazy 
cocycles and categorical structures, {\sl J. Pure Appl. Algebra} 
{\bf 209} (2007), 687--701.    

\bibitem{psvo}
F. Panaite, M. D. Staic, F. Van Oystaeyen, Pseudosymmetric braidings, 
twines and twisted algebras, {\sl J. Pure Appl. Algebra} 
{\bf 214} (2010), 867--884. 

\bibitem{ps}
F.  Panaite, M.  D. Staic, A quotient of the braid group related to pseudosymmetric braided categories, {\sl Pacific J. Math.} 
{\bf 244} (2010), 155--167.

\bibitem{pvo}
F. Panaite, F. Van Oystaeyen,  L-R-smash biproducts, double biproducts and 
a braided category of Yetter-Drinfeld-Long bimodules, 
{\sl Rocky Mount. J. Math.} {\bf 40} (2010), 2013--2024.

\bibitem{radford}
D. E. Radford, On Kauffmann's knot invariants arising from finite dimensional 
Hopf algebras, 
in ''Advances in Hopf algebras'', Lecture Notes Pure Appl. Math. Vol. 158, Dekker, 
New York, 1994, pp. 205--266.

\bibitem{doru}
M. D. Staic, Pure-braided Hopf algebras and knot invariants, 
{\sl J. Knot Theory Ramifications} {\bf 13} (2004), 385--400.


\end{thebibliography}
\end{document}